\documentclass[10pt,A4paper,oneside,reqno]{amsart}%

\usepackage{amssymb}
\usepackage{amsmath}%
\usepackage{amsfonts}%
\usepackage{graphicx}
\usepackage{color}

\newtheorem{theorem}{Theorem}[section]

\newtheorem{definition}[theorem]{Definition}
\newtheorem{example}[theorem]{Example}

\newtheorem{lemma}[theorem]{Lemma}

\newtheorem{proposition}[theorem]{Proposition}
\newtheorem{remark}[theorem]{Remark}

\numberwithin{equation}{section}

\def\R{\mathbb{R}}

\begin{document}

\title[Evolution equations with superlinear growth]{Evolution equations with nonlocal initial conditions and superlinear growth}



 \email{irene.benedetti@unipg.it} \email{simone.ciani@unifi.it}
 
  \maketitle

 \centerline{\scshape Irene Benedetti} \medskip {\footnotesize
 
   \centerline{Department of Mathematics and Computer Science}
   \centerline{University of Perugia, Italy}
 } 

 \medskip

 \centerline{\scshape Simone Ciani} \medskip {\footnotesize
  
   \centerline{Department of Mathematics and Informatics "Ulisse Dini"}
   \centerline{University of Florence, Italy} }

 \bigskip

\begin{abstract}

We carry out an analysis of the existence of solutions for a class of nonlinear partial differential equations of parabolic type. The equation is associated to a nonlocal initial condition, written in general form which includes, as particular cases, the Cauchy multipoint problem, the weighted mean value problem and the periodic problem. The dynamic is transformed into an abstract setting and by combining an approximation technique with the Leray-Schauder continuation principle, we prove global existence results. By the compactness of the semigroup generated by the linear operator, we do not assume any Lipschitzianity, nor compactness on the nonlinear term or on the nonlocal initial condition. In addition, the exploited approximation technique coupled to a Hartman-type inequality argument, allows to treat nonlinearities with superlinear growth. Moreover, regarding the periodic case, we are able to show the existence of at least one periodic solution on the half line.

\end{abstract}


\section{Introduction}

\noindent In this paper we consider the following class of nonlinear partial differential equations of parabolic type

\begin{equation}
\label{eq3}
u_t = \Delta u + h(t,x,u(t,x)) \quad \mbox{for} \quad (t,x) \in ]0,T[ \times \Omega
\end{equation}
where $\Omega \subset \mathbb{R}^k$, is a bounded domain with $C^2$-boundary and $h \colon [0,T] \times \Omega \times \mathbb{R} \to \mathbb{R}$ is a given map, under the following assumptions:
\begin{itemize}
\item[$(h1)$] for every $v \in \mathbb{R}$, $h(\cdot,\cdot,v)\colon [0,T] \times \Omega \to \mathbb{R}$ is a measurable function;
\item[$(h2)$] for every $t \in [0,T]$ and $x \in \Omega$, $h(t,x,\cdot) \colon \mathbb{R} \to \mathbb{R}$ is continuous;
\item[$(h3)$] there exists $m > 0$ and $\ell: [0,T] \times \Omega \to \mathbb{R}_+$ such that
\begin{itemize}
    \item[-] $\ell(\cdot,x) \in L^\infty([0,T];\mathbb{R}_+)$ for a.e. $x \in \Omega$;
    \item[-] $\ell(t,\cdot) \in L^q(\Omega;\mathbb{R}_+)$ for a.e. $t \in [0,T]$;
    \end{itemize} 
    \noindent and such that
$$ 
|h(t,x,v)| \leq \ell(t,x) + m |v|^{p/q}, \; \mbox{for every} \; v \in \mathbb{R} \; \mbox{and for a.e.} \; (t,x) \in [0,T] \times \Omega,
$$
with $2 \leq q < p < \infty$ for $k \leq 2$ and $2 \leq q < p < \infty$, $\displaystyle\frac{pq}{p-q} > \frac{k}{2}$, for $k > 2$;
\item[$(h4)$] $v \; h(t,x,v) \leq 0,$ for every $v \in \mathbb{R}$ and for a.e. $(t,x) \in [0,T] \times \Omega$.
\end{itemize}
\noindent The symbol $\Delta$ denotes the usual Laplace operator and we consider Dirichlet boundary conditions on $\partial\Omega$. 

\noindent A simple example of a map satisfying the required assumptions is 
\begin{equation}
\label{model example}
h(t,x,u)=-\frac{\sin(u)+2}{t^2+1} u^3,
\end{equation} 
see Example \ref{ex-finale}. More generally, $h$ can be a cubic polynomial, thus, the reaction diffusion model considered fits into the general Chafee-Infante class of partial differential equations studied by Henry \cite{H}. In such a model, the semilinear parabolic equation describes the evolution of the gene frequencies in the diffusion approximation for migration and selection at a multiallelic locus. Henry investigated the problem of existence of solutions as well as equilibria for two alleles in the case of homogeneous, isotropic migration (corresponding to the Laplacian). Since then, the problem was extended by many authors, for instance, Lou and Nagylaki in \cite{LN} considered the case of multiple alleles and of arbitrary migration (corresponding to an arbitrary elliptic operator); Huang and Huang in \cite{HH} consider \eqref{eq3} with $h(t,x,u)=\lambda(t)(u-u^3)$ and prove the existence of periodic mild solutions; Viorel in \cite{Vi} studies \eqref{eq3} with Neumann boundary conditions, combined to an integral-type nonlocal initial condition and proves the existence of global solutions near asymptotically stable equilibrium points.

\noindent We associate to the above equation a general nonlocal initial condition:
\begin{equation}
    \label{nonlocal-cond}
    u(0,\cdot)=g(u), 
 \end{equation} \noindent 
 where $g:C([0,T];L^p(\Omega; \mathbb{R})) \to L^p(\Omega; \mathbb{R})$, for $2 \leq p < \infty$, is a continuous function satisfying hypotheses $(g1)$, $(g2)$ of Section \ref{S2}. These assumptions allows to consider all the usual examples of nonlocal initial conditions such as periodic/antiperiodic, multipoint and integral-type initial conditions.

\noindent The growing interest for the existence of solutions of \eqref{eq3} which satisfy given nonlocal initial conditions arises with the possibility of these trajectories to capture additional information on the dynamics. In particular, the study of differential problems with nonlocal initial conditions was started by Byszewski with his pioneering paper \cite{B}, where the initial condition is given by
\begin{equation}
\label{g}
u(0)+g(t_1,\dots,t_m, u(\cdot))=u_0
\end{equation}
where $0 < t_1 < \dots < t_m \leq T$. Concrete nonlocal initial-boundary value problems for semilinear parabolic equations arising in physics (particularly in the mathematical modeling of heat conduction or diffusion processes) are analyzed in \cite{D}, \cite{J}. In particular, in \cite{D} the multipoint initial condition is used to describe the diffusion phenomenon of a small amount of gas in a transparent tube. In these cases, the condition \eqref{g} allows the measurements at $t=0, t_1,\dots, t_m$, rather than just at $t=0$. So more information is available. Subsequently, the abstract problem associated to \eqref{eq3} with initial nonlocal conditions has been extensively studied in the literature. For instance, Boucherif and Precup in \cite{BP} consider the multipoint condition under the assumption of the compactness of the semigroup generated by the linear part. In the same setting, see also \cite{ZCL} for functional semilinear differential equations and \cite{CLF} for strong solutions. Paicu and Vrabie \cite{PV} consider a general nonlocal initial condition of type \eqref{nonlocal-cond}.

\noindent Compared with the existing literature on the argument, the main novelty of our result is that we construct an abstract theory to solve differential equations that can take into account the superlinear growth of the Nemytskii operator associated to the function $h$, as well as the nonlocal initial condition \eqref{nonlocal-cond}. This technique was developed by the first author in \cite{BR} for an abstract Cauchy problem. In this paper we extend it to abstract nonlocal problems. More precisely, the approach in both \cite{BR} and here is constructed in a generic abstract framework, and its application is not limited to partial differential equations in the form \eqref{eq3}. 

\noindent In this abstract setting, to prove the existence of at least one solution, by exploiting the compactness of the semigroup generated by the Laplacian operator, we do not require that the nonlinearity be neither locally Lipschitz nor monotone, nor completely continuous, but only that it satisfies a Carath\'eodory condition.  

\noindent We remark that if $h(t,x,0) \neq 0$ for some $(t,x) \in [0,T] \times \Omega$, or $g(0) \neq 0$, then we obtain a non-zero solution of \eqref{eq3}.

\noindent Regarding the periodic case, we extend the obtained existence result to the half-line. More precisely, we prove the existence of at least one mild solution of the problem
\begin{equation} \label{periodic-infty}
\begin{cases}
u_t=\Delta u + h(t,x,u),& x \in \Omega,t>0,\\
u=0,& x \in \partial \Omega, t>0,\\
u(t,x)=u(t+T,x), & t \geq 0, \, x \in \bar{\Omega},
\end{cases}
\end{equation}
for a map $h: \R_+ \times \Omega \times \R \to \R$, $T$-periodic in the first variable and satisfying assumptions $(h1)$-$(h4)$ with $[0,T]$ replaced by $\R_+$, see Theorem \ref{ex-pde-infty}.

\noindent In relation to the above, in \cite{E}, \cite{GK}, \cite{Huska}, \cite{Quittner2} the authors prove existence results for positive classical solutions of \eqref{periodic-infty} under superlinear growth conditions on the nonlinear term, with different type of restrictions on the exponents $p$ and $q$, and in \cite{BPQ} for radial solutions. We also refer the reader to the monograph \cite{Quittner} dedicated to this topic. 

\noindent Finally, adding the following monotonicity type assumption
\begin{itemize}
\item[$(h5)$] $(u-v) \; (h(t,x,u)-h(t,x,v)) \leq 0$ for every $u,v \in \mathbb{R}$ and for a.e. $(t,x) \in [0,T] \times \Omega$,
\end{itemize}
we prove the uniqueness of the solution for the equation \eqref{eq3} associated to the nonlocal initial condition \eqref{nonlocal-cond}. Simple examples of a maps that satisfy conditions $(h5)$ and $(h4)$ are $-u^\alpha$ for each odd $\alpha \in \mathbb{N}$.

\section{Discussion and results}
\label{Discussion}

\noindent The usual approach to study the existence of solutions for the equation \eqref{eq3} consists in writing it as an abstract ordinary differential equation in a suitable infinite dimensional framework. More precisely, one considers the Nemytskii operator $f:[0,T] \times L^p(\Omega; \mathbb{R}) \to L^p(\Omega; \mathbb{R})$, associated to $h:[0,T] \times \Omega \times \mathbb{R} \to \mathbb{R}$ and the Laplace operator $A:D(A)=W^{1,p}_0(\Omega;\mathbb{R}) \cap W^{2,p}(\Omega;\mathbb{R}) \subset L^p(\Omega; \mathbb{R}) \to L^p(\Omega; \mathbb{R})$, $Aw=\Delta w$ for every $w \in D(A)$. Thus, one obtains the ordinary differential problem
\begin{equation}
\label{abstract-intro}
\begin{cases}
u'(t) = A u(t) + f(t,u(t)), \ \mbox{for a.e.}\ t\in [0,T], \\
u(0) = g(u).
\end{cases}
\end{equation}
The problem of finding mild solutions is then transformed into a suitable fixed point problem.

\noindent As stated in the well known Vainberg Theorem, the Nemytskii operator $f$ maps continuously the space $L^p(\Omega; \mathbb{R})$ into itself if and only if $h$ is sublinear.
\begin{theorem}[Theorem 19.1 in \cite{Vain}]
\label{th-Vain}
Let $B$ be a measurable set in a $s$-dimensional euclidean space and $h:B \times \mathbb{R} \to \mathbb{R}$ be a Carath\'eodory function, i.e. continuous with respect to $u \in \mathbb{R}$ for almost every $x \in B$ and measurable with respect to $x \in B$ for every $u \in \mathbb{R}$. Then the Nemytskii operator associated to $h$, $f(u)(x)=h(x,u(x))$, is a continuous and bounded operator from $L^p(B;\mathbb{R})$ into $L^q(B;\mathbb{R})$, $p,q \in [1,+\infty)$, if and only if there exist a function $a \in L^q(B;\mathbb{R})$ and a constant $b \geq 0$ such that for every $v \in \mathbb{R}$
$$
|h(x,v)| \leq a(x)+b |v|^{p/q}.
$$
\end{theorem}

\noindent Therefore, following the described abstract approach the nonlinearity term $h$ is forced to have a sublinear growth. 

\noindent  We overcome this difficulty by considering the domain and the arrival set of the Nemytskii operator $f$ two different spaces, namely $f:[0,T] \times L^p(\Omega; \mathbb{R}) \to L^q(\Omega; \mathbb{R})$, $2 \leq q < p < \infty$. Then, exploiting the fact that the semigroup generated by the Laplacian on $L^p(\Omega; \mathbb{R})$ can be extended to $L^q(\Omega; \mathbb{R})$ with $p > q$, by means of an approximation technique developed in \cite{XZF} (see also \cite{Reich}) and the Leray-Schauder continuation principle, we obtain the existence of a mild solution of \eqref{eq3} associated to \eqref{nonlocal-cond} localized in a ball of radius $R_0$ and center $0$. 

\noindent In case of lack of compactness of the semigroup generated by the linear part in \eqref{abstract-intro}, the usual assumption is to require that the nonlocal initial condition is compact, condensing or a contraction, e.g. see, among other papers, \cite{G-F}, \cite{X}, \cite{ZL}. Moreover, we mention also the papers \cite{HSA2} for $g$ depending also on the derivative of the solution $u$, \cite{H1}, \cite{V3} for delay equations and \cite{O},\cite{ZHL} for impulsive equations. The nonlocal problems studied under these assumptions do not recapture the periodic problem. In \cite{BLT} and in \cite{XCM}, the authors overcome this impediment by considering all the assumptions of regularity with respect to the weak topology. This approach, however, introduces restrictions on which forcing terms can be considered. In this paper, exploiting the compactness of the semigroup, we avoid these restrictions both on the nonlocal initial condition and on the forcing term. The idea is to require on the nonlinear term $f$ only a condition of Carath\'eodory type and on $g$ a condition slightly stronger than continuity, but weaker than compactness, see assumption $(g1)$ below. This idea derives from the observation that the nonlocal initial condition $g$ of type \eqref{g} is completely determined on $[\delta,T]$ for some small $\delta > 0$, i.e., such a $g$ ignores $t=0$. Liang Liu and Xiao in \cite{LLX}, see also \cite{CL}, \cite{HSA1} and \cite{XZF}, in the case of compact semigroup, generalize this idea and formulate a related condition for a general mapping $g:C([0,T];E) \to E$, with $E$ a Banach space:\vskip0.1cm \noindent 
\begin{itemize}
\item[$(g^*)$] $g$ is continuous and there is a $\delta \in (0,T)$ such that if $u,v \in C([0,T];E)$ are such that $u(s)\equiv v(s)$ for every $s \in [\delta,T]$ then $g(u)=g(v)$.
\end{itemize} \noindent Nevertheless, this hypothesis bears a strong requirement on the behaviour of $g$ within a nonzero measure interval around the starting time, that it is not satisfied by the classical mean value integral-type condition
$$
u(0,x)=\displaystyle \frac{1}{T} \int_0^T u(t,x) \, dt \quad \mbox{a.e. on} \; \Omega.
$$
In order to consider this nonlocal condition, or more general integral initial conditions, it is possible to require the following assumption:\vskip0.1cm \noindent 
\begin{itemize}
\item[$(g1)$] 
If $\{u_n\}_{n \in \mathbb{N}} \subset C([0,T];E)$ and $u_n(t) \rightarrow u(t),$ for $t \in (0,T]$, with $u \in C([0,T];E)$  then $g(u_n)\rightarrow g(u)$.
\end{itemize} \noindent To see that this assumption is weaker than the previous one, it is enough to consider a sequence  $\{u_n\}_{n \in \mathbb{N}}$ in $C([0,T];E)$ such that $u_n(t) \rightarrow u(t)$, for $t \in (0,T]$, with $u \in C([0,T];E)$, and consider the sequence $v_n:[0,T] \to E$ so defined:
\begin{equation*}
v_n(t)= \begin{cases}
u_n(\delta), & t\in [0,\delta],\\
u_n(t), & t \in (\delta,T].
\end{cases}  
\end{equation*}  \noindent
The sequence $\{v_n\}_{n \in \mathbb{N}}$ converges to $v:[0,T] \to E$ defined as 
\begin{equation*}
v(t)= \begin{cases}
u(\delta), & t\in [0,\delta],\\
u(t), & t \in (\delta,T].
\end{cases} 
\end{equation*}
So that, by assumption $(g^*)$, we obtain the desired condition $g(u_n)=g(v_n) \rightarrow g(v)=g(u)$. On the other hand, these conditions are not equivalent. For instance, the map $g$ which identifies the mean value condition, i.e.
$$
g(u)=\displaystyle\frac{1}{T}\int_0^T u(s) \,ds
$$
satisfies $(g1)$ but not $g^*$. To see it, consider $\delta \in (0,T)$ and $u:[0,T] \to \mathbb{R}$ defined as
$$
u(t)=\left\{\begin{array}{ll}
\displaystyle\frac{t}{\delta} & t \in [0,\delta) \\\\
1 & t \in [\delta,T]
\end{array}
\right.
$$
and $v:[0,T] \to \mathbb{R}$, $v \equiv 1$.

\noindent We refer the reader to \cite{BMT} and \cite{PV} for existence results for semilinear differential equations with nonlocal initial conditions satisfying $(g1)$ and to \cite{LY} for a discussion on the assumption $g^*$ and alternate condition both to $g^*$ and $(g1)$. Hypothesis $(g1)$ is general enough to consider a nonlocal initial condition belonging to one of the following map classes:
\begin{itemize}
\item[1.] {\it periodic/antiperiodic condition}: 
\begin{equation}
\label{periodic}
u(0,x)=u(T,x) \; / \; u(0,x)=-u(T,x), \quad \mbox{for every} \quad x \in \Omega;
\end{equation}
\item[2.] {\it multipoint condition} 
\begin{equation}
\label{multi}
u(0,x)=\displaystyle\sum_{i=1}^m c_i \gamma(u(t_i,x)), \quad \mbox{for every} \quad x \in \Omega,
\end{equation}
with $c_i \in \mathbb{R}$, $i=1,\dots,m$, $\gamma:\mathbb{R} \to \mathbb{R}$ and $0 < t_{1}<\cdots<t_{m}\leq T;$
\item[3.] {\it integral type condition}:
\begin{equation}
\label{integral}
u(0,x)=\displaystyle\int_{0}^T \eta(t,x,u(t,x))\,dt, \quad \mbox{for every} \quad x \in \Omega,
\end{equation}
where $\eta:[0,T] \times \Omega \times \mathbb{R} \to \mathbb{R}$.
\end{itemize}

\noindent Hence, we are able to prove the existence and uniqueness of mild solutions of \eqref{eq3}, according to Definition \ref{def-sol}, associated to each one of the above nonlocal initial conditions, as stated in the following theorems. 

\begin{theorem}
\label{ex-pde1}
Consider equation \eqref{eq3} under the assumptions $(h1)-(h4)$. Then the equation \eqref{eq3} associated to the nonlocal condition \eqref{periodic} admits at least one mild solution $u \in C([0,T];L^p(\Omega; \mathbb{R}))$ such that $\|u(t)\|_p < R_0$, for every $t \in [0,T]$ and for a suitable $R_0 > 0$.

\noindent Furthermore, if in addition we assume $(h5)$ and we require that 
\begin{equation}
\label{fixed-initial-condition}
u(0,x)=u_0(x)=u(T,x) \quad \mbox{for a.e.} \quad x \in \Omega
\end{equation}
with $u_0 \in L^p(\Omega; \mathbb{R})$, then the mild solution $u \in C([0,T];L^p(\Omega; \mathbb{R}))$ of \eqref{eq3}-\eqref{fixed-initial-condition} is unique.
\end{theorem}

\begin{theorem}
\label{ex-pde2}
Consider equation \eqref{eq3} under the assumptions $(h1)-(h4)$. 
Moreover, we assume that
\begin{itemize}
\item[$(i_1)$] $\displaystyle\sum_{i=1}^m c_i \leq 1$;
\item[$(i_2)$] $\gamma:\mathbb{R} \to \mathbb{R}$ is a continuous function such that $|\gamma(v)| \leq |v|$ for every $v \in \mathbb{R}$.
\end{itemize}
Then the equation \eqref{eq3} associated to the nonlocal condition \eqref{multi} admits at least one mild solution $u \in C([0,T];L^p(\Omega; \mathbb{R}))$ such that $\|u(t)\|_p < R_0$, for every $t \in [0,T]$ and for a suitable $R_0 > 0$.

\noindent Furthermore, if, instead of {\rm $(i_2)$}, we assume that
\begin{itemize}
\item[$(i_2')$] \begin{itemize}
\item[-] $\gamma:\mathbb{R} \to \mathbb{R}$ is Lipschitz continuous with a Lipschitz constant $\ell \in (0,1)$,
\item[-] $\gamma(0)=0$,
\end{itemize}
\end{itemize}
then the mild solution $u \in C([0,T];L^p(\Omega; \mathbb{R}))$ of \eqref{eq3}-\eqref{multi} is unique.
\end{theorem}

\begin{theorem}
\label{ex-pde3}
Consider equation \eqref{eq3} under the assumptions $(h1)-(h4)$. 
Moreover, we assume that
\begin{itemize}
\item[$({ii}_1)$] $\eta:[0,T] \times \Omega \times \mathbb{R} \to \mathbb{R}$ is a Carath\'eodory function 
\item[$({ii}_2)$] there exists $\alpha \in L^1([0,T];\mathbb{R})$ with $\|\alpha\|_{L^1([0,T])} \leq 1$ such that 
$$
|\eta(t,x,v)| \leq \alpha(t) |v| \quad \mbox{for a.e.} \quad (t,x) \in [0,T] \times \Omega \quad \mbox{and for every} \quad v \in \mathbb{R}.
$$
\end{itemize}
Then the equation \eqref{eq3} associated to the nonlocal condition \eqref{integral} admits at least one mild solution $u \in C([0,T];L^p(\Omega; \mathbb{R}))$ such that $\|u(t)\|_p < R_0$, for every $t \in [0,T]$ and for a suitable $R_0 > 0$.

\noindent Furthermore, if instead of {\rm $({ii}_2)$} we assume that
\begin{itemize}
\item[$({ii}_2')$] there exists $\alpha \in L^1([0,T];\mathbb{R})$ with $\|\alpha\|_{L^1([0,T])} < 1$ such that 
\begin{itemize}
    \item[-] $|\eta(t,x,v)-\eta(t,x,u)| \leq \alpha(t) |v-u|$ for a.e. $(t,x) \in [0,T] \times \Omega$ and for every $u,v \in \mathbb{R}$,
    \item[-] $\eta(t,x,0)=0$ for a.e. $(t,x) \in [0,T] \times \Omega$,
\end{itemize}
\end{itemize}
then the mild solution $u \in C([0,T];L^p(\Omega; \mathbb{R}))$ of \eqref{eq3}-\eqref{integral} is unique.

\end{theorem}

\begin{remark}
{\rm Notice that assumption $(h5)$ coupled with the hypothesis $h(t,x,0)=0$ for a.e. $(t,x) \in [0,T] \times x \in \Omega$, implies condition $(h4)$.}
\end{remark}

\noindent We point out the fact that, exploiting assumption $(g1)$, unlike the cited papers \cite{H1}, \cite{HSA2}, \cite{G-F}, \cite{O}, \cite{X}, \cite{ZHL}, \cite{ZL}, we are able to recover the periodic or the anti-periodic condition, as well as, unlike \cite{CL}, \cite{HSA1}, \cite{LLX}, and \cite{XZF} we can consider the mean value integral condition.

\vspace{.1cm}

\noindent Finally, we will prove the existence and uniqueness of at least one mild solution of \eqref{periodic-infty} $u \in C(\R_+, L^p(\Omega; \mathbb{R}))$, see Section \ref{sec-periodic}.

\noindent To the best of our knowledge, the most general result regarding the restrictions on the growth of $h$ for the problem \eqref{periodic-infty} is the one in \cite{Huska}, where the author assumes that $h(t,x,v)=m(t)\varphi(x,v)$ and
\begin{itemize}
\item[$(h3^*)$] $|\varphi(x,v)| \leq C(1+|v|^{\tilde{p}})$ with $1<\tilde{p}<(k+2)/(k-2)$ if $k>2$, and $1 < \tilde{p} < \infty$ if $k \leq 2$.
\end{itemize}
In comparison to this result, while on one side we find only mild solutions and do not take into account the sign problem, on the other side the growth condition $(h3)$ is less restrictive than $(h3^*)$. Moreover, our existence result is not limited to the periodic problem.

\vspace{.1cm}
\noindent 
The outline of this paper is the following. In Section \ref{S2} we present the abstract setting. In Sections \ref{App} and \ref{main-th} we prove an abstract existence result via the approximation technique. In Section \ref{existence-pde} we give the proof of our main results, i.e. the existence and uniqueness of mild solutions of \eqref{eq3} associated to \eqref{periodic}, \eqref{multi}, or \eqref{integral}. Finally in Section \ref{sec-periodic} we obtain the existence and uniqueness of at least one mild periodic solution on the half line.


\section{Abstract Theory} 

\label{S2}

\noindent As stated in the Introduction, we extend to differential problems with nonlocal initial conditions a theoretical approach developed in \cite{BR} for Cauchy differential problems in a generic abstract framework. For this reason, we consider the problem \eqref{abstract-intro} in the abstract setting determined by two generic real Banach spaces $(E, \|\cdot\|_E)$, $(F, \|\cdot\|_F)$ such that $E \subseteq F$ and assume that $E$ has a strictly convex dual $E^*$.

\noindent We denote with $B_{E,r}$ the ball in $E$ of radius $r>0$ and with $Q_r$ the ball in $C([0,T];E)$ of radius $r$ with the supremum norm, denoted by $\|\cdot\|_0$. 
For any $x,y \in E$, the normalized upper semi-inner product on $E$ 
$$
[u,v]_+=\lim_{h\downarrow 0} \, [u,v]_h,
$$
where for $h \in \mathbb{R}\setminus \{0\}$ we set
$$
[u,v]_h:=\displaystyle\frac{1}{h} \left (\|u+hv\|_E - \|u\|_E \right), 
$$
is well defined (see Lemma 1.4.1 and Definition 1.4.2 of \cite{V1}). Moreover, denoting by $\langle\cdot,\cdot\rangle$ the duality product between $E^*$ and $E$ and by $J_E:E \multimap E^*$ the duality map, i.e.
\begin{equation}\label{JDef}
J_E(u)=\{u^* \in E^* \colon \|u^*\|_{E^*}=\|u\|_E \; \mbox{and} \; \langle u^*, u \rangle = \|u\|_E^2 \}
\end{equation}
for every $u,v \in E,\, u \neq 0$, we have
$$
[u,v]_+ = \displaystyle\frac{1}{\|u\|_E} \sup\{\langle u^*,v\rangle \colon u^* \in J_E(u)\},
$$
see Lemma 1.4.2 and 1.4.3 in \cite{V1}. 

\noindent In particular, since $E^*$ is strictly convex, $J$ is a single valued map. Thus, we get
\begin{equation}
\label{inner-product}
[u,v]_+ = \displaystyle\frac{1}{\|u\|}_{E} \langle J_E(u),v \rangle.
\end{equation}

\noindent Moreover, for every $\lambda \in \mathbb{R}$ and $x \in X$, it holds
$$
J_E(\lambda u)=\lambda J_E(u).
$$

\noindent We consider the abstract differential problem \eqref{abstract-intro} under the following assumptions: \medskip

\noindent $A\colon D(A)\subset E\to E$ is a linear operator such that

\begin{itemize}
	
\item[$(A1)$] $0 \in D(A)$, generating a compact $C_0-$semigroup of contractions $\{S(t)\}_{t\geq 0}$ in $E$;
\item[$(A2)$] the semigroup $\{S(t)\}_{t\geq 0}$ can be extended to a semigroup in $F$, i.e. 
	\begin{itemize}
	\item[(A2.i)] there exists a semigroup $\{S^*(t)\}_{t\geq 0}$ on $F$ generated by $A$ such that for every $w \in E$, it holds $S^*(t)w=S(t)w$;
	\item[(A2.ii)] for every $v \in F$ and $t > 0$, we have $S^* (t)v \in E$;
	\item[(A2.iii)] there exists a function $c \in L^r([0,T];\mathbb{R}_+)$, with $1 \leq r \leq \infty$ such that for any $v \in F$ it holds 
	\begin{equation*}
	\label{norm}
	\|S^*(t) v\|_E \leq c(t) \|v\|_F \; \mbox{for every} \; t \in (0,T];
	\end{equation*}
	\end{itemize}

\end{itemize} 	\medskip
	
the map $f:[0,T] \times E \to F$ is specified as follows
	
\begin{itemize}
	
	\item[$(f1)$] for every $v \in E$ the map $f(\cdot,v):[0,T]  \to F$ is measurable;
	
	\item[$(f2)$] for a.e. $t \in [0,T]$ the map $f(t,\cdot):E \to F$ is continuous;
		
	\item[$(f3)$] for every bounded subset $D\subset E$ there exists a function $\nu_{D}\in L^{r^\prime}([0,T];\mathbb{R}_+)$, with $\frac{1}{r}+\frac{1}{r^\prime}=1$ and $r^\prime=\infty$ if $r =1$, such that
								$$
									{\|f(t,v)\|}_{F}\leq \nu_{D}(t),
								$$
								for a.e. $t\in [0,T]$ and all $v\in D$;
								
	\item[$(f4)$] there exist constants $0 < r_0 < R_0$ and $n_0 \in \mathbb{N}$ such that for every $n > n_0$
                           $$
                           {\langle J_E(v),S^*\left(\frac{1}{n}\right) f(t,v) \rangle} \leq 0,
                           $$
                           for a.e.  $t \in [0,T]$ and for every $v\in E$ such that $r_0 < {\|v\|}_{E} < R_0$, where $J_E$ denotes the duality map on $E$;	
\end{itemize}  \medskip

\noindent and finally $g:C([0,T];E) \to E$ satisfies the conditions

\begin{itemize}

\item[$(g1)$] if $\{u_n\}_{n \in \mathbb{N}} \subset C([0,T];E)$ and $u_n(t) \rightarrow u(t)$, $t \in (0,T]$, with $u \in C([0,T];E)$ for $t \in (0,T]$ then $g(u_n)\rightarrow g(u)$;
\item[$(g2)$] $\displaystyle\sup_{u \in Q_{R}}\| g(u) \|_E \leq {R}$, for $r_0 < R < R_0$, where $r_0,R_0>0$ are given by the transversality condition $(f4)$ on $f$ above.

\end{itemize}

\begin{remark}
{\rm Because of condition $(A2)$, we can denote the $C_0$-semigroup generated by $A$, on the space $E$ or on the space $F$, by the very same symbol $\{S(t)\}_{t\geq 0}$.

\noindent Moreover notice that, being $\{S(t)\}_{t\geq 0}$ a $C_0$-semigroup on the space $F$, there exists a constant $M > 0$ such that
\begin{equation}
\label{M}
\|S(t)\|_F \leq M,
\end{equation}
for every $t \in [0,T]$.
 }
\end{remark}

\noindent We will prove the existence of at least one mild solution of \eqref{abstract-intro}, i.e. a function $u$ satisfying the following definition.

\smallskip

\begin{definition}
\label{def-sol}
By a solution of \eqref{abstract-intro} we mean a function $u \in C([0,T]; E)$ satisfying the nonlocal initial condition \eqref{nonlocal-cond} and such that for each $t \in [0,T]$
\begin{equation} \label{mildsolution}
u(t)=S(t)g(u) + \int_0^t S(t-\tau) f(\tau,u(\tau)) \,  d\tau \;.
\end{equation} 
\end{definition} 
\noindent Now we consider the ball
\begin{equation} \label{Qr}
Q_R= \{ q \in C([0,T];E): ||q(t)||_{E} \leq R, \forall t \in [0,T] \}, \quad \quad R \in (r_0,R_0),
\end{equation} 
where $r_0$ and $R_0$ are from assumption $(f4)$. By {\it a solution on $Q_R$} we mean a map $u \in Q_R$ satisfying \eqref{mildsolution} for each $t \in [0,T]$.

\noindent Moreover to prove uniqueness of the solution we assume that
\begin{itemize}
\item[$(f5)$] there exists $n_0 \in \mathbb{N}$ such that for every $n > n_0$
                           $$
                           {\langle J_E(u-v),S^*\left(\frac{1}{n}\right)(f(t,u)-f(t,v)) \rangle} \leq 0,
                           $$
                           for a.e.  $t \in [0,T]$ and for every $u,v\in E$, where $J_E$ denotes the duality map on $E$;	
\item[$(g2')$] there exists $L \in (0,1)$ such that for every $u,v\in E$
$$
\|g(u)-g(v)\|_E \leq L \|u-v\|_0
$$
and
$$
L R + \|g(0)\|_E \leq R \quad \mbox{for} \quad r_0 < R < R_0,
$$
where $r_0,R_0>0$ are given by the transversality condition $(f4)$ on $f$ above.
\end{itemize} 

\noindent Our main abstract result is the following.

\begin{theorem} \label{main}
Let conditions $(A1),(A2),(f1)-(f4)-(g1)-(g2)$ being satisfied, then the set of solutions on $Q_R$ is nonempty and compact in $C([0,T];E)$.

\noindent Moreover, if condition $(f5)$ is satisfied and condition $(g2)$ is replaced by $(g2')$ the solution is unique.
\end{theorem}

\noindent There are several definitions of solutions of \eqref{abstract-intro}. In the proof of our main abstract existence result we make use of the concept of integral solution and of the equivalence between the integral and the mild solutions in the particular case of linear Cauchy problems. More precisely, consider the linear problem 
\begin{equation}\label{eq0}
\begin{cases}
u'(t) = A u(t) + \beta(t), \ \mbox{for a.e.}\ t\in [0,T], \\
u(0) = \xi \in E,
\end{cases}
\end{equation}
where $A\colon D(A)\subset E\to E$ is the infinitesimal generator of a semigroup $\{S(t)\}_{t \geq 0}$ and $\beta \in L^1([0,T];E)$ is a given map. 
\begin{definition}[Definition 3.5.1 in \cite{LL} and Definition 1.7.4 in \cite{V1}]
\label{sol-integral}
A function $u:[0,T] \to E$ is called an integral solution of \eqref{eq0} on $[0,T]$ if $u \in C([0,T];E)$ satisfies $u(0)=\xi$ and 
$$
\|u(t)-x\|_E \leq \|u(s)-x\|_E + \displaystyle\int_s^t [u(\tau)-x,\beta(\tau)+Ax]_+ \, d\tau
$$
for each $x \in D(A)$ and $0 \leq s \leq t \leq T$.
\end{definition}
\noindent From Theorems 1.7.3, 1.7.4 and 1.8.2 in \cite{V1} we can deduce the following existence result and the equivalence between the mild and the integral solutions.
\begin{theorem}
\label{existence-limit}
For each $\xi \in E$ and $\beta \in L^1([0,T];E)$ there exists a unique integral solution $u$ of \eqref{eq0} on $[0,T]$.
\end{theorem}
\begin{theorem}
\label{mild-integral}
A function $u:[0,T] \to E$ is a mild solution of \eqref{eq0} if and only if $u$ is an integral solution of \eqref{eq0} on $[0,T]$ satisfying $u(0)=\xi$.
\end{theorem}

\noindent The following estimate holds for integral solutions of \eqref{eq0} and it will be useful to prove the uniqueness of the mild solutions.
\begin{theorem}[Theorem 1.7.5 in \cite{V1}]
\label{estimate}
Let $\beta_1,\beta_2 \in L^1([0,T];E)$ and let $u,v \in C([0,T];E)$ two solutions of the equation in \eqref{eq0} corresponding to $\beta_1$ and to $\beta_2$ respectively. Then
$$
\|u(t)-v(t)\|_E \leq \|u(s)-v(s)\|_E + \displaystyle\int_s^t [u(\tau)-v(\tau),\beta_1(\tau)-\beta_2(\tau) ]_+ \, d\tau
$$
for each $0 \leq s \leq t \leq T$.
\end{theorem}

\noindent Moreover, the proof of the existence result is based on an approximation technique and on the compactness result Proposition \ref{compactness} (see \cite{BR}).

\noindent Given $\xi \in E$ and $\beta \in L^{r^\prime}([0,T];F)$, where $1 \leq r^\prime < \infty$ is such that $\frac{1}{r}+\frac{1}{r^\prime}=1$, where $1 \leq r \leq \infty$ is defined in (A2), and $r^\prime=\infty$ if $r =1$, we denote by $\mathcal{F}(\xi,\beta):[0,T] \to E$ the mild solution of \eqref{eq0}, that is
\begin{equation}
\label{mild}
	\mathcal{F}(\xi,\beta)(t) = S(t)\xi + \int_0^t S(t-s) \beta(s)\, ds, \ \mbox{for every}\ t\in [0,T].
\end{equation}
\begin{proposition}[Proposition 3.4 in \cite{BR}]
\label{compactness}
If $A:D(A) \subset E \to E$ satisfies $(A1)$ and $(A2)$, then for each bounded subset $B$ of $E$ and each subset $H$ in $L^{r^\prime}([0,T];F)$ such that $\{\|\beta\|_F^{r^\prime}, \beta \in H\}$ is uniformly integrable, the set $\mathcal{F}(B \times H)$ is relatively compact in $C([\delta,T];E)$ for each $\delta \in (0,T)$. If, in addition, $B$ is relatively compact in $E$, then $\mathcal{F}(B \times H)$ is relatively compact in $C([0,T];E)$.
\end{proposition}
\noindent We transform the problem to find solutions of \eqref{abstract-intro} into a fixed point problem and we apply the Leray-Schauder continuation principle, see e.g. \cite{FP} or the original paper \cite{LS}.
\begin{theorem}
\label{cp}
Let $Q$ be a closed subset of a Banach space $\mathcal{B}$ and let $\Sigma \colon Q\times [0,1]\to \mathcal{B}$ be a continuous map sending bounded subsets of $Q \times [0,1]$ into relatively compact subsets of $\mathcal{B}$. Assume that
\begin{itemize}
\item[(a)] $\Sigma(x,0)=x_0 \in \mbox{int}(Q), \quad \forall \; x \in Q$;
\item[(b)] The fixed point set 
               $$F=\{x \in Q, x=\Sigma(x,\lambda), \; \mbox{for some} \; \lambda \in [0,1] \}$$
               is bounded and does not meet the boundary $\partial Q$ of Q.
\end{itemize}
Then the map $x \mapsto \Sigma(x,1)$ has a fixed point in $Q$.
\end{theorem}               

\section{Approximating problem}
\label{App}

\noindent The proof of Theorem \ref{main} relies on an approximation technique. We introduce a family of approximating problems: for $n \in \mathbb{N}$, we consider the following semilinear problem.
\begin{equation} \label{Pn}
(P_{n})\begin{cases}
u'(t)=A u(t)+ S(1/n) f(t,u(t)),& \text{for a.e.}\quad  t \in [0,T],\\
u(0)=S(1/n)g(u) \in E
\end{cases}
\end{equation} 
\begin{lemma} \label{approx}
Let conditions $(A1),(A2),(f1)-(f4)-(g1)-(g2)$ being satisfied, then there exists $n_0 \in \mathbb{N}$ such that for every $R \in (r_0,R_0)$ and every $n>n_0$ the problem $P_n$ admits a mild solution $u_n \in C([0,T];E)$ satisfying $ || u_n(t)||_E \leq R,$ for every $t \in [0,T]$.
\end{lemma}
\begin{proof}
Our strategy is to show that for every $n>n_0$ it is possible to apply the Leray-Schauder continuation principle to each problem $(P_n)$. 

\noindent Let $n_0$ from assumption $(f4)$ and let $n > n_0$ be fixed. We define the operator $\Sigma_n: Q_R \times [0,1] \rightarrow C([0,T];E)$ by
\begin{equation} \label{sigma}
\Sigma_{n} (q,\lambda)(t)= \lambda S(t) S(1/n) g(q)+ \lambda \int_0^t S(t-\tau) S(1/n) f(\tau,q(\tau)) d\tau \quad t \in [0,T].
\end{equation}\noindent
The operator $\Sigma_n$ is well defined because of condition (A2.ii)  (see for instance Proposition 3.3 in \cite{BR}), and a mild solution of problem $(P_n)$ is exactly a fixed point of the operator $\Sigma_n(\cdot,1)$. 

\vspace{.1cm}
\noindent As a first step we show that the operator $\Sigma_n$ is continuous.

\noindent Let $q \leftarrow \{ q_k \}_{k \in \mathbb{N} }\subset Q_R$ and $\lambda \leftarrow \{ \lambda_k \}_{k \in \mathbb{N} } \subset [0,1] $ be two convergent sequences. We observe that by assumptions $(g1)$ and $(g2)$ on the map $g$ and by the properties of the semigroup $\{S(t)\}_{t \geq 0}$ it follows that
$$
\begin{array}{lcl}
|| \lambda_k S(t) S(1/n) g(q_k)- \lambda S(t) S(1/n) g(q)||_E & \leq & || \lambda_k S(t) S(1/n) (g(q_k)-g(q))||_E \\
&& +||(\lambda_k-\lambda) S(t)S(1/n)g(q)||_E \\
& \leq & ||g(q_k)-g(q)||_E+ |\lambda_k-\lambda| \|g(q)\|_E \\
& \leq & ||g(q_k)-g(q)||_E+ R |\lambda_k-\lambda| \rightarrow 0, \quad \text{for} \, k \, \rightarrow \infty.
\end{array}
$$
Next, by $(f2)$ we have that
$$
\|f(t,q_k(t)) - f(t,q(t))\|_F \to 0 \quad \forall \; t \in [0,T],
$$
hence, by (A2.iii) it follows
$$
\|S(t)S\left(\frac{1}{n}\right)(f(t,q_k(t))-f(t,q(t)))\|_E \leq c\left(\frac{1}{n}\right) \|f(t,q_k(t)) - f(t,q(t))\|_F \to 0 \quad \forall \; t \in [0,T].
$$
Moreover, by (A2.iii) and $(f3)$ we get
$$
\|S(t)S\left(\frac{1}{n}\right) f(t,q_k(t))\|_E \leq c\left(\frac{1}{n}\right) \|f(t,q_k(t))\|_F \leq c\left(\frac{1}{n}\right) \nu_{B_R}(t) \quad \mbox{for a.e.} \; t \in [0,T].
$$
Thus, by the Lebesgue's Dominated Convergence Theorem we conclude that for every $t \in [0,T]$
\begin{align*}
\|\Sigma_n(q_k,\lambda_k)(t)-\Sigma_n(q,\lambda)(t)\|_E & \leq \displaystyle || \lambda_k S(t) S(1/n) g(q_k)- \lambda S(t) S(1/n) g(q)||_E +\\
& \ \ \ \ \displaystyle + |\lambda_k - \lambda| \int_0^t \left \|S(t-\tau) S\left(\frac{1}{n}\right) f(\tau,q(\tau))\right \|_E \,d\tau +\\
& \ \ \ \ \displaystyle + \lambda_k \int_0^t \left \|S(t-\tau) S\left(\frac{1}{n}\right) (f(\tau,q_k(\tau))-f(\tau,q(\tau)))\right\|_E \,d\tau \\
& \ \ \ \ \leq \displaystyle ||g(q_k)-g(q)||_E + |\lambda_k-\lambda| R+\displaystyle |\lambda_k - \lambda| c\left(\frac{1}{n}\right) T^{\frac{1}{r}}\|\nu_{B_R}\|_{L^{r^\prime}([0,T],\mathbb{R}_+)}+\\
& \ \ \ \ \displaystyle + \lambda_k \int_0^T \left \|S(t-\tau) S\left(\frac{1}{n}\right) (f(\tau,q_k(\tau))-f(\tau,q(\tau)))\right \|_E \,d\tau \to_{k \to \infty} 0
\end{align*}                                                                                            
Hence $\Sigma_n(q_k,\lambda_k) \to \Sigma_n(q,\lambda)$ in $C([0,T];E)$, obtaining the continuity of the operator $\Sigma_n$.

\vspace{.1cm}
\noindent Now, as a second step, we show that for every $n \in \mathbb{N}$, the operator $\Sigma_n$ sends $Q_R \times [0,1]$ into a relatively compact set of $C([0,T];E)$.\\
\noindent We observe that $\Sigma_n(Q_R \times [0,1])(0)$ is a relatively compact set, since it coincides with $[0,1] \times S(1/n)g(Q_R)$ and $g(Q_R)$ is bounded by $(g2)$, while $S(1/n)$ is a compact operator. 

\noindent Finally, by $(f3)$ and the boundedness of the semigroup $\{S(t)\}_{t \geq 0}$ in the Banach space $F$, there exists a function $\nu_{B_R} \in L^{r^\prime}([0,T];\mathbb{R}_+)$ such that 
$$
\left\|S \left(\frac{1}{n}\right) f\left(t,q(t)\right)\right\|_F \leq M \nu_{B_R}(t),\quad \mbox{for a.e.} \; t \in [0,T] \; \mbox{and for every}\, q \in Q_R,
$$
where $M$ is defined in \eqref{M}, implying that the set $\left\{S\left(\frac{1}{n}\right) f\left(\cdot,q(\cdot)\right), q \in Q_R \right\}$ is a family of maps in $L^{r^\prime}([0,T],F)$ such that $\left\{\left \|S\left(\frac{1}{n}\right) f\left(\cdot,q(\cdot)\right) \right\|_F^{r^\prime},\, q \in Q_R \right\}$ is uniformly integrable. Therefore, observing that $\Sigma_n(q,\lambda)=\lambda \mathcal{F}(S(1/n)g(q),S\left(\frac{1}{n}\right)f((\cdot),q(\cdot)))$ for every $(q,\lambda)\in Q_R \times [0,1]$ we obtain, by Proposition \ref{compactness}, that the set $\Sigma_n(Q_R \times [0,1])$ is relatively compact in $C([0,T];E)$. 

\vspace{.2cm}

\noindent To show that $\Sigma_n(Q_R \times \{0\}) \subset \mbox{int}(Q_R)$, it is enough to observe that $\Sigma_n(Q_R \times \{0\})\equiv 0$.

\vspace{.2cm}

\noindent Finally we need to prove that the operator $\Sigma_n(\cdot,\lambda)$ has no fixed points on $\partial Q_R$ for every $\lambda \in [0,1]$ and $n > n_0$, where $n_0$ is from $(f4)$.\newline
We argue by contradiction: let us assume that there exists $\overline{\lambda} \in [0,1]$, $\overline{u} \in Q_R$ and $t_0 \in [0,T]$ such that $\overline{u}=\Sigma_n(\overline{u},\overline{\lambda})$ and $\|\overline{u}(t_0)\|_E = R$. Since $\overline{\lambda}=0$ implies $\overline{u} \equiv 0$ and $\overline{\lambda}=1$, gives the existence of at least one fixed point $\overline{u}=\Sigma_n(\overline{u},1)$, we may assume $\overline{\lambda}\in (0,1)$. 

\noindent Notice that $t_0 \neq 0$. Indeed, if $t_0=0$ we have by hypothesis $(g2)$
$$
R=\|\overline{u}(0)\|_E=\|\Sigma_n(\overline{u},\overline{\lambda})(0)\|_E = \overline{\lambda} \| S(1/n) g(\bar{u})\|_E < R.
$$
Hence, there exists $\delta > 0$ such that $r_0 < \|\overline{u}(t)\|_E \leq R$ for every $t \in [t_0 -\delta,t_0]$ and $\|\overline{u}(t_0-\delta)\|_E < R$. 

\vspace{.1cm}

\noindent Denoting by $f_n(t)=S\left(\frac{1}{n}\right)f(t,\overline{u}(t))$, $g_n=S(1/n)g(\overline{u})$, $t \in [0,T]$, we consider the linear problem
\begin{equation}\label{eq2}
\begin{cases}
u'(t) = A u(t) + f_n(t), \ \mbox{for a.e.}\ t\in [0,T], \\
u(0) = g_n \in E.
\end{cases}
\end{equation}
By the fact that $\|\overline{u}(t)\|_E \leq R$ for every $t \in [0,T]$, by $(A2)$ and $(f3)$, we have that 
$$
\|f_n(t)\|_E =\left\|S\left(\frac{1}{n}\right) f\left(t,\overline{u}(t)\right)\right\|_E \leq c\left(\frac{1}{n}\right) \nu_{B_R}(t) \quad \mbox{for a.e.} \; t \in [0,T],
$$
obtaining that $f_n \in L^{r^\prime}([0,T];E)$. Let $u \in C([0,T];E)$ be the unique mild solution of \eqref{eq2}, i.e.
$$
u(t)=S(t) g_n+ \displaystyle \int_0^t S(t-s) f_n(s) \,ds, \quad t \in [0,T].
$$
By Theorem \ref{mild-integral}, we have that $u$ is the unique integral solution of \eqref{eq2}, i.e.
$$
\|u(t)-x\|_E \leq \|u(s)-x\|_E + \displaystyle \int_s^t [u(\tau)-x, f_n(\tau)+Ax]_+ \, d\tau
$$
for each $x \in D(A)$ and $0 \leq s \leq t \leq T$. Since $E$ has a strictly convex dual, bearing in mind that 
\[ [x,y]_+= \frac{1}{||x||}_{E} \left \langle J_E(x),y \right\rangle, \]
we have that
$$
\|u(t)-x\|_E \leq \|u(s)-x\|_E + \displaystyle \int_s^t \frac{1}{\|u(\tau)-x\|_E} \langle J_E(u(\tau)-x), f_n(\tau)+Ax \rangle\, d\tau
$$
for each $x \in D(A)$ and $0 \leq s \leq t \leq T$. By the definition of the operator $\Sigma_n(\cdot,\overline{\lambda})$ and the fact that $\overline{u}$ is a fixed point of it, for every $t \in [0,T]$, we obtain
$$
\begin{array}{lcl}
u(t)=\displaystyle S(t)g_n + \displaystyle \int_0^t S(t-\tau) f_n(\tau) \,ds & = & \displaystyle S(t)S\left(\frac{1}{n}\right) g(\overline{u}) + \displaystyle \int_0^t S(t-\tau) S\left(\frac{1}{n}\right)f\left(\tau,\overline{u}(\tau)\right) \,d\tau = \displaystyle\frac{\overline{u}(t)}{\overline{\lambda}}.                                                                                
\end{array}
$$
Now, considering $x=0 \in D(A)$ and observing that $\|u(s)\|_E > 0$ for every $s \in [t_0-\delta,t_0]$, it follows that
$$
\begin{array}{lcl}
0 & < & \displaystyle\frac{\|\overline{u}(t_0)\|_E - \|\overline{u}(t_0-\delta)\|_E}{\overline{\lambda}} \vspace{.1cm}\\
   & = & \|u(t_0)\|_E - \|u(t_0-\delta)\|_E \leq \displaystyle \int_{t_0 - \delta}^{t_0} \frac{1}{\|u(\tau)\|_E} \langle J_E(u(\tau)),f_n(\tau) \rangle \,d\tau \vspace{.1cm}\\
   & = & \displaystyle \int_{t_0 - \delta}^{t_0} \frac{1}{\|u(\tau)\|_E} \left\langle J_E(u(\tau)),S\left(\frac{1}{n}\right) f\left(\tau,\overline{u}(\tau)\right) \right\rangle \, d\tau \vspace{.1cm}\\
   & = & \displaystyle \int_{t_0 - \delta}^{t_0} \frac{1}{\|u(\tau)\|_E} \left\langle J_E\left(\frac{\overline{u}(\tau)}{\overline{\lambda}}\right),S\left(\frac{1}{n}\right) f\left(\tau,\overline{u}(\tau)\right) \right \rangle \, d\tau \vspace{.1cm}\\
   & = & \displaystyle \int_{t_0 - \delta}^{t_0} \frac{1}{\overline{\lambda}\, \|u(\tau)\|_E} \left\langle J_E(\overline{u}(\tau)),S\left(\frac{1}{n}\right) f\left(\tau,\overline{u}(\tau)\right) \right\rangle\,d\tau,
\end{array}
$$
where the last equality follows from the properties of the duality map. Now, by $(f4)$ for every $n > n_0$ we get the contradiction
$$
0 <  \displaystyle\frac{\|\overline{u}(t_0)\|_E - \|\overline{u}(t_0-\delta)\|_E}{\overline{\lambda}} \leq 0.
$$
\noindent By Leray-Schauder continuation principle, for every $n > n_0$, we obtain the existence of a fixed point $u = \Sigma_n(u,1)$. Thus for every  $n > n_0$, we get a mild solution of $(P_n)$.   
\end{proof}
\begin{lemma} \label{unique-approx}
Let conditions $(A1),(A2),(f1)-(f5)-(g1)-(g2')$ being satisfied, then there exists $n_0 \in \mathbb{N}$ such that for every $R \in (r_0,R_0)$ and every $n>n_0$ the problem $P_n$ admits a unique mild solution $u_n \in C([0,T];E)$ satisfying $ || u_n(t)||_E \leq R,$ for every $t \in [0,T]$.
\end{lemma}
\begin{proof}
First of all, notice that condition $(g2')$ implies condition $(g2)$. Indeed
$$
\displaystyle\sup_{u \in Q_{R}} \| g(u) \|_E \leq \displaystyle\sup_{u \in Q_{R}} \left (\|g(u)-g(0)\|_E + \|g(0)\|_E\right) \leq \displaystyle\sup_{u \in Q_{R}} \left( L \|u\|_E + \|g(0)\|_E\right) \leq L R +\|g(0)\|_E \leq R.
$$
Thus, existence follows from Lemma \ref{approx}. Now, assume by contradiction that there exist two mild solutions $u, v \in C([0,T];E)$ of $(P_n)$ with $u \neq v$. Hence, in particular, $u$ and $v$ are integral solutions of the equation in \eqref{eq0} corresponding to $\beta_1, \beta_2 \in L^1([0,T],E)$ defined as $\beta_1(t)=S\left(\frac{1}{n}\right)f(t,u(t))$ and $\beta_2(t)=S\left(\frac{1}{n}\right)f(t,v(t))$, $t \in [0,T]$ respectively. Let
$$
t_0=\inf\left\{t \in [0,T] \; \mbox{such that} \; u(t) \neq v(t) \right\}.
$$
By continuity, $t_0 \neq T$ and there exists $\delta > 0$ such that $u(t) \neq v(t)$ for every $t \in [t_0, t_0+\delta]$. 

\noindent Firstly, assume $t_0 \neq 0$, thus, $u(t)=v(t)$ for every $t \in [0,t_0)$ and $u(t) \neq v(t)$ for every $t \in [t_0,t_0+\delta]$. Let $t \in (t_0,t_0+\delta]$, by Theorem \ref{estimate} and $(f5)$, we have
$$
\begin{array}{lcl}
0 < \|u(t)-v(t)\|_E & \leq & \|u(0)-v(0)\|_E + \displaystyle\int_0^t [u(\tau)-v(\tau),\beta_1(\tau)-\beta_2(\tau)]_+\, d\tau \vspace{.1cm}\\ 
& = & \displaystyle\int_{t_0}^t [u(\tau)-v(\tau),\beta_1(\tau)-\beta_2(\tau)]_+\, d\tau \vspace{.1cm}\\ 
& = & \displaystyle\int_{t_0}^t \frac{1}{\|u(\tau)-v(\tau)\|_E} \langle J_E(u(\tau)-v(\tau)),S\left(\frac{1}{n}\right)(f(\tau,u(\tau))-f(\tau,v(\tau))) \rangle \, d\tau \leq 0
\end{array}
$$
getting a contradiction.
In the case $t_0=0$, for every $t \in [0,\delta]$, by Theorem \ref{estimate}, $(f5)$ and $(g2')$ we have
$$
\begin{array}{lcl}
0 < \|u(t)-v(t)\|_E & \leq & \|u(0)-v(0)\|_E + \displaystyle\int_0^t [u(\tau)-v(\tau)),\beta_1(\tau)-\beta_2(\tau)]_+\, d\tau \vspace{.1cm}\\ 
& = & \|g(u)-g(v)\|_E \vspace{.1cm}\\
&& + \displaystyle\int_{0}^t \frac{1}{\|u(\tau)-v(\tau)\|_E} \langle J_E(u(\tau)-v(\tau)),S\left(\frac{1}{n}\right)(f(\tau,u(\tau))-f(\tau,v(\tau))) \rangle \, d\tau \vspace{.1cm}\\
& \leq & L \|u-v\|_0 \vspace{.1cm}\\
&& + \displaystyle\int_{0}^t \frac{1}{\|u(\tau)-v(\tau)\|_E} \langle J_E(u(\tau)-v(\tau)),S\left(\frac{1}{n}\right)(f(\tau,u(\tau))-f(\tau,v(\tau))) \rangle \, d\tau \vspace{.1cm}\\
& \leq & L \|u-v\|_0.
\end{array}
$$
Thus we have
$$
0 < \displaystyle\sup_{t \in [0,T]} \|u(t)-v(t)\| =\|u-v\|_0 \leq L \|u-v\|_0 < \|u-v\|_0,
$$
getting again a contradiction. As a consequence, we have for every  $n > n_0$, the uniqueness of the solution of $P_n$.
\end{proof}
\section{Proof of Theorem \ref{main}}
\label{main-th}
\noindent By Lemma \ref{approx}, we know that there exists $n_0 > 0$ such that problems $(P_n)$ have at least one mild solution for every $n > n_0$. We consider now the set of these mild solutions. More precisely, by the characterization introduced in Lemma \ref{approx}, we consider the set 
$$
M_R=\{u_n \in C([0,T];E)\cap Q_R \; : \; u_n=\Sigma_n(u_n,1), n > n_0\}.
$$
Let $u_n \in M_R$. By the fact that $u_n \in Q_R$, we have that $\|u_n(t)\|_E \leq R$ for every $t \in [0,T]$. Thus, by $(f3)$ there exist a function $\nu_{B_R} \in L^{r^\prime}([0,T];\mathbb{R}_+)$ 
$$
\left\|S \left(\frac{1}{n}\right) f\left(t, u_n(t)\right)\right\|_F \leq M \nu_{B_R}(t),\quad \mbox{for a.e.} \; t \in [0,T],
$$
implying that the set $G=\left\{S\left(\frac{1}{n}\right) f\left(\cdot,u_n(\cdot)\right), n > n_0 \right\}$ is a family of maps in $L^{r^\prime}([0,T],F)$ such that $\left\{\left \|S\left(\frac{1}{n}\right) f\left(\cdot,u_n(\cdot)\right) \right\|_F^{r^\prime},\, n > n_0 \right\}$ is uniformly integrable. Hence, applying Proposition \ref{compactness} with the bounded set $B=\{S(1/n) g(M_R),\, n > n_0\}$ and $G$,  for each chosen $\delta>0$ we have the relative compactness of $M_R$ in $C([\delta,T];E)$. 

\noindent  Let $\{ u_n \}_{n \in \mathbb{N}}$ be a sequence in $M_R$. By previous considerations, for each $\delta>0$ chosen, there exists $\{u_n^{\delta}\}_{n \in \mathbb{N}} \subset \{u_n \}_{n \in \mathbb{N}} \rightarrow u_{\delta}^*$ in $C([\delta,T];E)$. By uniqueness of the limit, when $\delta_{n}<\delta_{n-1}$ then $u^*_{\delta_{n}}= u^*_{\delta_{n-1}}$ in $[\delta_{n-1},T]$. Consider any sequence of numbers $\delta_n \downarrow 0$ and perform a Cantor diagonal argument to show that there exists a function $u^* \in C((0,T];E)$ such that a particular (diagonal) sub-sequence $\{ u_n^n\}_{n \in \mathbb{N}}$ of each $\{u_n^{\delta_n}\}_{n \in \mathbb{N}}$ converges to $u^*$ in $C((0,T];E)$. Moreover, since $E$ is reflexive and $\{ u_n \}_{n \in \mathbb{N}}\subset Q_R$ by property $(g2)$ there exists $u_0 \in E$ such that $g(u_n^n) \rightharpoonup u_0$. Now we define a function $\bar{u} \in C([0,T];E)$ by
\[
\bar{u}(t):= S(t)u_0+ \int_0^t S(t-s) f(s,u^*(s)) \, ds.
\] 
\noindent We claim that $\{u_n^n(t) \}_{n \in \mathbb{N}} \rightharpoonup \bar{u}(t)$ for each $t \in [0,T]$.
By the continuity of $S\left(\frac{1}{n}\right)$ and of $f(t, \cdot)$, for every $n \in \mathbb{N}$ we obtain
$$
S\left(\frac{1}{n}\right) f\left(t,u_n^n(t)\right) \stackrel{F}\to f(t,u^*(t)), \quad \mbox{for a.e.} \; t \in (0,T],
$$
moreover, since $\|u_n^n(t)\|_E \leq R$ for every $t \in [0,T]$, the convergence is dominated
$$
\|S\left(\frac{1}{n}\right) f\left(t,u_n^n(t)\right)\|_F \leq M \nu_{B_R}(t) \quad \mbox{for a.e.} \; t \in [0,T],
$$
Thus we get
$$
 \displaystyle\int_0^t S(t-\tau) S\left(\frac{1}{n}\right) f\left(\tau,u_n(\tau)\right) \, d\tau \to \displaystyle\int_0^t S(t-\tau) f(\tau,u^*(\tau))\, d\tau, \quad \forall \; t \in [0,T].
$$
Again, by the continuity and linearity of the semigroup we have that
$$
S(t)S\left(\displaystyle\frac{1}{n}\right)g(u_n^n)  \rightharpoonup S(t)S\left(\displaystyle\frac{1}{n}\right)u_0
$$
and so the claimed result. By the uniqueness of the limit we have that $u^*(t)=\bar{u}(t)$ for every $t \in (0,T]$, thus $\{u_n^n(t) \}_{n \in \mathbb{N}} \rightarrow \bar{u}(t)$ for each $t \in (0,T]$, with $\bar{u} \in C([0,T];E)$. By property $(g1)$ this means that $g(u_n^n) \rightarrow g(\bar{u})$ in $C([0,T];E)$, declaring again by uniqueness that $g(\bar{u})=u_0$. In conclusion we get that for every $t \in [0,T]$
$$
\bar{u}(t):= S(t)g(\bar{u})+ \int_0^t S(t-s) f(s,\bar{u}(s)) \, ds,
$$
proving that $\bar{u}$ is a solution on $Q_R$ of \eqref{abstract-intro}.

\noindent So the set of solutions on $Q_R$ 
\[\mathcal{S}:=\mathcal{F}(g(\mathcal{S}),f([0,T],\mathcal{S}))\] \noindent
is nonempty, now we prove that it is compact by the very same argument used above for $M_R$. Indeed, $f([0,T],\mathcal{S})$ is uniformly integrable by $(f3)$ and $g(\mathcal{S})$ is bounded by $g(2)$, so, by Proposition \ref{compactness}, for each $\delta>0$ we obtain that $\mathcal{F}((g(\mathcal{S}),f(s,\mathcal{S})))$ is relatively compact in $C([\delta,T];E)$. Then again by a Cantor diagonal argument we have that for each sequence $\{u_n \}_{n \in \mathbb{N}} \subseteq \mathcal{S}$ there exists a subsequence $\{u_{n}^n\}_{n \in \mathbb{N}}$ converging to $u^* \in C((0,T];E)$ and exactly as before we can continuously extend $u^*$ to $[0,T]$. Next, by $(g1)$ we have that $g(u_n^n) \to g(u^*)$ and so we get that $g(\mathcal{S})$ is relatively compact. Applying again Proposition \ref{compactness} we obtain that $\mathcal{S}$ is a relatively compact set.

\noindent Assuming, in addition, condition $(h5)$ and $(g2')$ instead of $(g2)$, by Lemma \ref{unique-approx} we have the existence of a unique sequence $\{u_n\}$, $n > n_0$, of solutions of the problems $(P_n)$. By the above reasonings we have that $\{u_n\}$, $n > n_0$ converges to some $\overline{u}$ solution on $Q_R$ of \eqref{abstract-intro}, getting the uniqueness of the solution. 


\section{Existence and uniqueness results for the problem \eqref{eq3}-\eqref{nonlocal-cond}} 
\label{existence-pde}

\noindent In this Section we prove the main result of the paper, i.e. the existence of at least one solution $u \in C([0,T];L^p(\Omega; \mathbb{R}))$ for the problem \eqref{eq3} associated to each one of the nonlocal initial conditions \eqref{periodic}, \eqref{multi}, \eqref{integral}. \vskip0.2cm \noindent
To prove such existence results we will apply Theorem \ref{main}. Indeed, the equation \eqref{eq3} associated to one of the conditions \eqref{periodic}, \eqref{multi}, \eqref{integral}, can be rewritten as an abstract evolution equation of the form \eqref{abstract-intro} with $E=L^p(\Omega;\mathbb{R})$ and $F=L^q(\Omega;\mathbb{R})$, $2 \leq q < p < \infty$. 

\noindent The Laplace operator $A:D(A) \subset L^p(\Omega;\mathbb{R}) \to L^p(\Omega;\mathbb{R})$ subjected to Dirichlet boundary conditions on $L^p(\Omega; \mathbb{R})$ and defined by 
$$
\begin{array}{l}
D(A)=W^{1,p}_0(\Omega;\mathbb{R}) \cap W^{2,p}(\Omega;\mathbb{R}), \\
A w=\Delta w,
\end{array}
$$ 
satisfies conditions $(A1)$ and $(A2)$.

\noindent Indeed, $A$ is the generator of a $C_0$-semigroup of contractions $\{S_p(t)\}_{t \geq 0}$ (see e.g. Theorem 4.1.3 and Remark 4.1.2 of \cite{V2}). Moreover, by Lemma 7.2.1 of \cite{V2}, for each $p,q \in [1,+\infty]$, each $\xi \in C(\overline{\Omega};\mathbb{R})$ and each $t \geq 0$, we have $S_p(t)\xi=S_q(t)\xi.$ Thus, we can denote the $C_0-$semigroup generated by the Laplace operator  subjected to the Dirichlet boundary conditions on any of the spaces $L^p(\Omega;\mathbb{R})$ by the very same symbol $\{S(t)\}_{t \geq 0}$. By Theorem 7.2.5 of \cite{V2}, $\{S(t)\}_{t \geq 0}$ is a compact semigroup. Finally, by Theorem 7.2.6 of \cite{V2}, for each $1 \leq q \leq p \leq \infty$, each $\xi \in L^q(\Omega;\mathbb{R})$, and each $t > 0$, we have 
$$
\|S(t)\xi\|_p \leq (4\pi t)^{-\frac{k}{2}\left(\frac{1}{q}-\frac{1}{p}\right)} \|\xi\|_q.
$$
Hence, $A$ satisfies $(A1)$ and $(A2)$ with $c(t)=(4\pi t)^{-\frac{k}{2}\left(\frac{1}{q}-\frac{1}{p}\right)}$. Notice that $\frac{k}{2}\left(\frac{1}{q}-\frac{1}{p}\right) < 1$, provided $2 \leq q < p < \infty$ for $k \leq 2$ and $2 \leq q < p < \infty$, $\displaystyle\frac{pq}{p-q} > \frac{k}{2}$, for $k > 2$, hence the function $c \in L^1([0,T],\mathbb{R}_+)$. \vspace{.2cm}

\subsection{Proof of Theorem \ref{ex-pde1}}

\noindent We will prove that all the hypotheses of Theorem \ref{main} are satisfied.

\noindent By $(h1),(h2),(h3)$ and the Vainberg Theorem (see \cite{Vain}) we have that the Nemytskii operator $f:[0,T] \times L^p(\Omega; \mathbb{R}) \to L^q(\Omega; \mathbb{R})$ defined as $f(t,u)(x)=h(t,x,u(x))$ maps the space $L^p(\Omega; \mathbb{R})$ into $L^q(\Omega; \mathbb{R})$ and is continuous. Moreover again by $(h3)$, we get
$$
\|f(t,u)\|_q = \displaystyle \bigg( \int_\Omega |h(t,x,u(x))|^q \, dx\bigg)^{\frac{1}{q}} \leq C (\|\ell(t,\cdot)\|_q^q + m \|u\|_p^p),
$$
by Minkovski inequality, where $C > 0$ is a suitable constant. Hence for every bounded subset $D$ of $L^p(\Omega; \mathbb{R})$, we have that 
$$
\|f(t,u)\|_q \leq C (\|\ell(t,\cdot)\|_q^q + C_1):=\nu_D(t),
$$
for a.e. $t \in [0,T]$ and for every $u \in D$, with $C_1 > 0$ another suitable constant. So assumption $(f3)$ is satisfied with $\nu_D \in L^\infty([0,T];\mathbb{R}_+)$.

\noindent Now, let $0 \ne u \in L^p(\Omega; \mathbb{R})$ and $ t \in [0,T]$. Denoting by $\text{sign}(x)$ the sign function and with
\[
[u>0]= \{x \in \Omega | u(x)>0  \},
\]
\[
[u<0]= \{x \in \Omega | u(x)<0  \},
\]
we perform the calculation to show that the transversality property $(f4)$ holds. First of all we recall the definition of the duality map in the space $L^p(\Omega; \mathbb{R})$. For every $0 \ne u \in L^p(\Omega; \mathbb{R}) $, we have
$$
\langle J_{L^p(\Omega; \mathbb{R})}(u), v \rangle = \displaystyle \frac{1}{\|u\|^{p-2}_p} \int_\Omega |u(\xi)|^{p-2} u(\xi) v(\xi) \,d\xi,
$$
see e.g. Example 1.4.4 in \cite{V1}. Let $N_0 \subset [0,T]$ and $\Omega_0 \subset \Omega$ be two sets with Lebesgue measure zero and let $(t,x) \in \{[0,T] \setminus N_0\} \times \{\Omega \setminus \Omega_0\}$ be such that $(h4)$ is satisfied. Thus, for every $u \in L^p(\Omega; \mathbb{R})$, $\|u\|_p \neq 0$ we have
\begin{equation*}
    \begin{aligned}
    \langle J_{L^p(\Omega; \mathbb{R})} (u),& S(1/n) f(t,u)\rangle=\\
    & \frac{1}{||u||_p^{p-2}} \displaystyle \int_{\Omega}
 |u|^{p-2}u \, S\left (\frac{1}{n}\right)(h(t,\xi,u(\xi)))\,  d\xi \\
 & = \displaystyle\frac{1}{||u||^{p-2}_p}  \bigg( \int_{[u>0]} |u|^{p-1} \, S\left (\frac{1}{n}\right) h(t,\xi,u(\xi)) \, d\xi -  \int_{[u<0]} |u|^{p-1}  \, S\left (\frac{1}{n}\right) h(t,\xi,u(\xi)) \, d\xi \bigg)\leq 0.
\end{aligned}
\end{equation*} 
\noindent We explain the last inequality. The semigroup $\{ S(t) \}$ generated by the Laplace operator subjected to the Dirichlet boundary conditions is a positive semigroup in $L^p(\Omega; \mathbb{R})$ when $p\ge2$, i.e. for every $w \in L^p(\Omega; \mathbb{R})$ such that $w(x) \geq 0$ for a.a. $x \in \Omega$
$$
S(t)w(x) \geq 0 \; \mbox{for a.e.} \; x \in \Omega \; \mbox{and} \; \mbox{for every} \; t \in [0,T],
$$
see \cite{V2}, Lemma 7.2.3 and moreover by property $(h4)$ we have
\[
h(t,\xi,u(\xi)) \leq 0, \quad \text{for a.e.} \quad \xi \in [u>0]
\] and the converse inequality in $[u<0]$.   

\noindent Finally, if the equation \eqref{eq3} is associated to a periodic or anti periodic condition, trivially satisfies conditions $(g1),(g2)$. 

\noindent Thus all the assumptions of Theorem \ref{main} are satisfied and we get the existence of at least one mild solution of  \eqref{eq3}-\eqref{periodic}.

\noindent Now we show that condition $(h5)$ on $h$ implies condition $(f5)$ on the superposition operator $f$. Indeed, similarly as above, for every $u, v \in L^p(\Omega; \mathbb{R})$, $\|u-v\|_p \neq 0$ we have
\begin{equation*}
    \begin{aligned}
    \langle J_{L^p(\Omega; \mathbb{R})} (u-v),& S(1/n) (f(t,u)-f(t,v))\rangle=\\
    & \frac{1}{||u-v||_p^{p-2}} \displaystyle \int_{\Omega}
 |u-v|^{p-2}(u-v) \, S\left (\frac{1}{n}\right)(h(t,\xi,u(\xi))-h(t,\xi,u(\xi)))\,  d\xi \\
 & = \displaystyle\frac{1}{||u-v||^{p-2}_p}  \bigg( \int_{[(u-v)>0]} |u-v|^{p-1} \, S\left (\frac{1}{n}\right) ((h(t,\xi,u(\xi))-h(t,\xi,u(\xi))) \, d\xi \\
 & -  \int_{[(u-v)<0]} |u-v|^{p-1}  \, S\left (\frac{1}{n}\right) (h(t,\xi,u(\xi))-h(t,\xi,u(\xi))) \, d\xi \bigg)\leq 0.
\end{aligned}
\end{equation*} 
Finally, since by \eqref{fixed-initial-condition}, $u(0,x)=v(0,x),\,x \in \Omega$ for every $u,v$ solutions of \eqref{eq3}, we get uniqueness of mild solutions of \eqref{abstract-intro} under conditions $(f1)-(f5)$ and $(g1)$, $(g2)$ and so we get the uniqueness of mild solution of \eqref{eq3}-\eqref{fixed-initial-condition}. Thus the proof is completed.

\subsection{Proof of Theorem \ref{ex-pde2}}

The only difference with the previous theorem is the nonlocal initial condition. Thus we have only to prove that conditions $(g1)$, $(g2)$ are satisfied.

\noindent In this case  the initial datum $g:C([0,T];L^p(\Omega; \mathbb{R})) \to L^p(\Omega; \mathbb{R})$ is defined as 
$$
g(v)(x)=\displaystyle\sum_{i=1}^m c_i \gamma(v(t_i)(x)).
$$
Clearly it is a continuous map and we have by $(i_1)$, $(i_2)$
$$
\displaystyle\sup_{v \in Q_{R}} \|g(v)\|_p \leq \sup_{v \in Q_{R}} \sum_{i=1}^m c_i \|v(t_i)\|_p \leq R,
$$
for $r_0 < R < R_0$. Moreover, if we consider $\{u_n\} \subset C([0,T];E)$ such that $u_n(t) \to u(t)$ for every $t \in (0,T]$, with $u \in C([0,T];E)$ we have that 
$$
\|g(u_n)-g(u)\|_p^p \leq \displaystyle\sum_{i=1}^m c_i \displaystyle\int_\Omega |\gamma(u_n(t_i)(x))-\gamma(u(t_i)(x))|^p \,dx.
$$
Now, as $u \in C([0,T];L^p(\Omega; \mathbb{R}))$ the sequence $u_n$ is equibounded for $n$ big and using (i) joint to the Dominated Convergence Theorem we obtain that $\gamma(t,x,u_n(t)(x))$ converges to $\gamma(t,x,u(t)(x))$ in $L^p(\Omega; \mathbb{R})$ for a.e $t \in [0,T]$ and so that $g(u_n) \to g(u)$ in $L^p(\Omega; \mathbb{R})$ as well.

\noindent Furthermore, if we replace $(i_2)$ by $(i_2')$ we get
$$
\begin{array}{lcl}
\|g(u)-g(v)\|_p & \leq & \displaystyle\sum_{i=1}^m c_i \left(\displaystyle\int_\Omega |\gamma(u(t_i)(x))-\gamma(v(t_i)(x))|^p \,dx\right)^{1/p} \\
&\leq& \displaystyle\sum_{i=1}^m c_i \left(\displaystyle\int_\Omega \ell^p |(u(t_i)(x))-(v(t_i)(x))|^p \,dx\right)^{1/p} \leq \ell \displaystyle\sup_{t \in [0,T]} \|u(t)-v(t)\|_p.
\end{array}
$$
So, $g:C([0,T];L^p(\Omega; \mathbb{R})) \to L^p(\Omega; \mathbb{R})$ is Lipschitz continuous with Lipschitz constant $L= \ell$. Finally, $g(0)=0$ and so trivially $LR+\|g(0)\|_p \leq R$, thus $(g2')$ is satisfied.

\noindent Hence, also in this case all the assumptions of Theorem \ref{main} are satisfied and we have the existence and uniqueness of mild solution for the problem \eqref{eq3}-\eqref{multi} on $Q_R$.

\subsection{Proof of Theorem \ref{ex-pde3}}

\noindent In this case the initial datum $g:C([0,T];L^p(\Omega; \mathbb{R})) \to L^p(\Omega; \mathbb{R})$ is written as 
$$
g(v)(x)=\displaystyle \int_{0}^T \eta(t,x,v(t)(x)) \, dt.
$$
By condition $({ii}_1)$ and $({ii}_2)$ we have 
$$
\begin{array}{lcl}
\displaystyle\sup_{v \in Q_{R}} \|g(v)\|_p^p & \leq & \displaystyle \sup_{v \in Q_{R}} \int_{\Omega} \bigg| \int_{0}^T \alpha(t) |v(t)(x)| \,dt \bigg|^p \, dx \vspace{.1cm}\\
& \leq & \displaystyle \sup_{v \in Q_{R}} \int_{\Omega} \bigg( \max_{[0,T]} |v(t)(x)| \bigg)^p\, \bigg| \int_{0}^T \alpha(t) \,dt \bigg|^p \, dx \vspace{.1cm}\\
& \leq & \displaystyle R^p \|\alpha\|_{L^1([0,T])}^p  \leq R^p.
\end{array}
$$
for $r_0 < R < R_0$. 

\noindent To verify condition $(g1)$, consider a sequence $v_n \in C([0,T];L^p(\Omega; \mathbb{R}))$ converging to a function $v^*\in C([0,T];L^p(\Omega; \mathbb{R}))$ for each $t>0$. Then, as $v^* \in C([0,T];L^p(\Omega; \mathbb{R}))$ the sequence $v_n$ is equibounded for $n$ big and using (i') joint to the Dominated Convergence Theorem we obtain that $\eta(t,x,v_n(t)(x))$ converges to $\eta(t,x,v^*(t)(x))$ in $L^p(\Omega; \mathbb{R})$ for a.e $t \in [0,T]$. Thus, we get
$$
\begin{array}{lcl}
\|g(v_n)-g(v^*)\|_p^p & = & \displaystyle\int_{\Omega} \left |\int_0^T (\eta(t,x,v_n(t)(x))-\eta(t,x,v^*(t)(x)))\,dt \right |^p \, dx \\
& \leq & c \displaystyle\int_0^T \left(\int_\Omega \left |\eta(t,x,v_n(t)(x))-\eta(t,x,v^*(t)(x)))\right |^p \, dx\right)\,dt 
\end{array}
$$
where the last term tends to zero by the convergence of $\eta(t,x,v_n(t)(x))$ to $\eta(t,x,v^*(t)(x))$ in $L^p(\Omega; \mathbb{R})$.

\noindent Furthermore, if we replace $({ii}_2)$ with $({ii}_2')$ we get
$$
\begin{array}{lcl}
\displaystyle \|g(u)-g(v)\|_p & \leq & \displaystyle\left(\int_{\Omega} \bigg| \int_{0}^T \alpha(t) |u(t)(x)-v(t)(x)| \,dt \bigg|^p \, dx\right)^{1/p} \vspace{.1cm}\\
& \leq & \displaystyle \left(\int_{\Omega} \bigg( \max_{[0,T]} |u(t)(x)-v(t)(x)| \bigg)^p\, \bigg| \int_{0}^T \alpha(t) \,dt \bigg|^p \, dx\right)^{1/p} \vspace{.1cm}\\
& \leq & \|\alpha\|_{L^1([0,T])}  \displaystyle\sup_{t \in [0,T]} \|u(t)-v(t)\|_p.
\end{array}
$$
So, $g:C([0,T];L^p(\Omega; \mathbb{R})) \to L^p(\Omega; \mathbb{R})$ is Lipschitz continuous with Lipschitz constant $L=\|\alpha\|_{L^1([0,T])}$. Finally, since $\eta(t,x,0)=0$ for a.e. $(t,x) \in [0,T] \times \Omega$ we get condition $(g2')$ as well.

\noindent Thus, also in this case all the assumptions of Theorem \ref{main} are satisfied and we have the existence and uniqueness of mild solution for the problem \eqref{eq3}-\eqref{integral}.

\vspace{.2cm}
\noindent We conclude this study with an example that include the mean value initial condition.

\begin{example}
\label{ex-finale}
{\rm An example of nonlocal differential problem that satisfies all the requirements is the following.
$$
\left\{\begin{array}{l}
u_t=\Delta u-\displaystyle\frac{\sin(u)+2}{1+t^2} u^3, \quad (t,x) \in (0,T) \times \Omega, \\
u(t,x)=0 \quad \mbox{a.e. on} \; (0,T) \times \partial \Omega, \\
u(0,x)=\displaystyle \frac{1}{T} \int_0^T \alpha(t) |u(t,x)| \, dt \quad \mbox{a.e. on} \; \Omega,
\end{array}
\right.
$$
where $\Omega \subset \mathbb{R}^k$, $k \geq 2$, is as in \eqref{eq3}, $\alpha \in L^1([0,T],\mathbb{R}_+)$ is such that $\|\alpha\|_1 \leq 1$. Indeed trivially, $h(t,x,u)=-\frac{\sin(u)+2}{1+t^2} u^3$ is a continuous function, assumptions $(h2)$ is satisfied for instance for $\frac{k}{3} < q < \infty$, $p=3q$, with $\ell \equiv 0$ and $m=3$. Finally,
$$
u h(t,x,u)=-\frac{\sin(u)+2}{1+t^2} u^4 \leq 0 \quad \mbox{for every} \quad u \in \mathbb{R}.
$$
\noindent 
So, applying Theorem \ref{ex-pde3} we obtain the existence of at least one global mild solution $u \in C([0,T];L^p(\Omega; \mathbb{R}))$.
}

\end{example}

\begin{example}
\label{ex-finale-2}
{\rm If we consider a similar equation as in Example \ref{ex-finale} associated to \eqref{fixed-initial-condition} with $u_0 \neq 0$ we get the uniqueness of a non zero mild solution. Consider the problem:
$$
\left\{\begin{array}{l}
u_t=\Delta u-u^3, \quad (t,x) \in (0,T) \times \Omega, \\
u(t,x)=0 \quad \mbox{a.e. on} \; (0,T) \times \partial \Omega, \\
u(0,x)=u_0(x)=u(T,x) \quad \mbox{a.e. on} \; \Omega,
\end{array}
\right.
$$
where $\Omega \subset \mathbb{R}^k$, $k \geq 2$, is as in \eqref{eq3}. As in Example \ref{ex-finale} it is possible to prove that conditions $(h1)-(h4)$ are satisfied. Moreover for every $u,v \in \mathbb{R}$
$$
(u-v)(-u^3+v^3)=(u-v)(v-u)(u^2+ uv + v^2)=-(u-v)^2 (u^2+ uv + v^2) \leq 0,
$$
so, condition $(h5)$ is satisfied as well and we get the existence and uniqueness of a global mild solution $u \in C([0,T];L^p(\Omega; \mathbb{R}))$ by Theorem \ref{ex-pde1}.}
\end{example}

\section{Periodic solutions on the half line}
\label{sec-periodic}

\noindent In this Section we prove an existence result for solutions of \eqref{periodic-infty}. To this aim, we have to consider the abstract problem \eqref{abstract-intro} in the half line. Let $T > 0$ be fixed and let $f:[0,+\infty) \times L^p(\Omega; \mathbb{R}) \to L^q(\Omega; \mathbb{R})$, $2 \leq q < p < \infty$, be the Nemytskii operator associated to a $T$-periodic function $h:[0,+\infty) \times \Omega \times \mathbb{R} \to \mathbb{R}$, thus obtaining the ordinary differential periodic problem
\begin{equation}
\label{Abstract2}
\begin{cases}
u'(t) = A u(t) + f(t,u(t)), \ \mbox{for a.e.}\ t\in \R_+, \\
u(t+T) = u(t), \quad \forall \, t \in \R_+
\end{cases}
\end{equation}
where $A$ is as before the Laplace operator. We consider the abstract differential problem \eqref{Abstract2} under the same assumptions of previous problem \eqref{abstract-intro}, with the exception of the following needed adjustments: 

\begin{itemize}
	\item[$(A2\infty)$] the semigroup $\{S(t)\}_{t\geq 0}$ can be extended to a semigroup in $F$, i.e. 
	\begin{itemize}
	\item[(A2.i)] there exists a semigroup $\{S^*(t)\}_{t\geq 0}$ on $F$ generated by $A$ such that for every $w \in E$, it holds $S^*(t)w=S(t)w$;
	\item[(A2.ii)] for every $v \in F$ and $t > 0$, we have $S^* (t)v \in E$;
	\item[(A2.iii$\infty$)] there exists a function $c \in L^r_{{\tiny\mbox{loc}}}(\R_+;\mathbb{R}_+)$, with $1 \leq r \leq \infty$ such that for any $v \in F$ it holds 
	\begin{equation*}
	\label{normS}
	\|S^*(t) v\|_E \leq c(t) \|v\|_F \; \mbox{for every} \; t \in (0,+\infty);
	\end{equation*}
\end{itemize}	
	\item[$(f1\infty)$] the map $f: \R_+\times E  \to F$ is Carath\'eodory;
	
	\item[$(f2\infty)$] for every $z \in E$ the map $f(\cdot,z):[0,+\infty) \to E$ is $T$-periodic;
	
	\item[$(f3\infty)$] for every bounded subset $D\subset E$ there exists a function $\nu_{D}\in L^{r^\prime}_{\tiny{\mbox{loc}}}(\R_+;\mathbb{R}_+)$, with $\frac{1}{r}+\frac{1}{r^\prime}=1$ and $r^\prime=\infty$ if $r =1$, such that
								\[
									{\|f(t,v)\|}_{F}\leq \nu_{D}(t),\quad \text{for a.e.} \quad t\in  \R_+, \quad \forall v\in D;
								\]
								
	\item[$(f4\infty)$] there exist constants $0 < r_0 < R_0$ and $n_0 \in \mathbb{N}$ such that for every $n > n_0$
                           $$
                           {\langle J_E(v),S^*\left(\frac{1}{n}\right) f(t,v) \rangle} \leq 0,
                           $$
                           for a.e.  $t \in [0,T]$ and for every $v\in E$ such that $r_0 < {\|v\|}_{E} < R_0$.
         \item[$(f5\infty)$] there exists $n_0 \in \mathbb{N}$ such that for everyfor every $n > n_0$
                           $$
                           {\langle J_E(v),S^*\left(\frac{1}{n}\right) (f(t,u)-f(t,v)) \rangle} \leq 0,
                           $$
                           for a.e.  $t \in [0,T]$ and for every $u,v\in E$.
                                 
 \end{itemize}  \medskip
\begin{remark}
\noindent {\rm Notice that by the construction of the solution on the half line (see the proof of Theorem \ref{main2}), we need only to assume $(f4\infty)$ for a.e. $t \in [0,T]$ and not, as to be expected, for a.e. $t \in \R_+$.}
\end{remark}                     
\noindent We will prove the existence of at least one mild solution of \eqref{Abstract2}, i.e. a periodic function $u \in C(\R_+, L^p(\Omega; \mathbb{R}))$ such that for each $t \in \R_+$ 
it holds \eqref{mildsolution}.
\begin{theorem}
 \label{main2}
Let conditions $(A1)$-$(A2\infty)$ and $(f1\infty)$-$(f4\infty)$ being satisfied, then there exists at least one solution $u \in C(\R_+;E)$ on $Q_R$ of the periodic problem \eqref{Abstract2}. 

\noindent Moreover if we assume $(f5\infty)$ and
\begin{equation}
\label{fixed-initial-condition-abstract}
u(0)=u_0=u(T)
\end{equation}
with $u_0 \in E$, then the solution $u \in C([0,T];E)$ of \eqref{Abstract2}-\eqref{fixed-initial-condition-abstract} is unique.

\end{theorem}

\begin{proof} By Theorem \ref{main}, there exists at least one mild solution of \eqref{abstract-intro}, $u:[0,T] \to E$, such that $u(0)=u(T)$. Now we show how to extend $u$ to $\R_+$ in order to obtain a solution of \eqref{Abstract2}. For every $t \in \R_+ = \cup_{m \in \mathbb{N}, m \geq 1} [(m-1)T,mT]$, there exists $m \in \mathbb{N}$ such that $t \in [(m-1)T,mT]$, so we define the map $\overline{u}:\R_+ \to E$ as
$$
\overline{u}(t)=u(t-(m-1)T),\quad m \geq 1,\; t \in [(m-1)T,mT]
$$
and we will prove that $\overline{u}$ is a periodic map satisfying \eqref{mildsolution} for every $t \in \R_+$.

\noindent Let $t \in  [(m-1)T,mT]$, then $t+T \in  [(mT,(m+1)T]$ and so
$$
\overline{u}(t+T)=u(t+T-mT)=u(t-(m-1)T)=\overline{u}(t),
$$
thus obtaining the periodicity of $\overline{u}$.
Now we proceed by an induction process. For $m=1$ is trivial, since $\overline{u}\equiv u$ in $[0,T]$ and $u$ is a mild solution of \eqref{abstract-intro}. Now, we assume that $\overline{u}$ satisfies \eqref{mildsolution} for every $t \in [(m-2)T,(m-1)T]$ and we will prove that this is still the case for every $t \in [(m-1)T,mT]$. Let $t \in \mathbb{R}_+$, and $m \in \mathbb{N}$ such that $t\in [(m-1)T,mT]$. Then $t=r+T$ with $r \in [(m-2)T,(m-1)T]$. So, by the periodicity of $\overline{u}$ and the inductive assumption we have that
$$
\begin{array}{lcl}
\overline{u}(t) & = & \overline{u}(r+T)=\overline{u}(r)=S(r)\overline{u}(T)+\displaystyle\int_0^r S(r-s) f(s,\overline{u}(s))\, ds \\
& = & S(r)\left[S(T)\overline{u}(T)+\displaystyle\int_0^T S(T-s) f(s,\overline{u}(s))\, ds\right] + \displaystyle\int_0^r S(r-s) f(s,\overline{u}(s))\, ds \\
& = & S(r+T)\overline{u}(T)+\displaystyle\int_0^T S(r+T-s) f(s,\overline{u}(s))\, ds+\displaystyle\int_0^r S(r-s) f(s+T,\overline{u}(s+T))\, ds \\
& = & S(r+T)\overline{u}(T)+\displaystyle\int_0^T S(r+T-s) f(s,\overline{u}(s))\, ds+\displaystyle\int_T^{r+T} S(r+T-\eta) f(\eta,\overline{u}(\eta))\, d\eta \\
& = & S(t)\overline{u}(T)+\displaystyle\int_0^t S(t-s) f(s,\overline{u}(s))\, ds.
\end{array}
$$

\noindent Hence $\overline{u}$ is a solution of problem \eqref{Abstract2}. 

\noindent Any two mild solutions $u, v \in C([0,T];E)$ of  problem \eqref{abstract-intro} with $g:C([0,T];E) \to E$ given by $g(u)=u(T)$, satisfying also \eqref{fixed-initial-condition-abstract}, are such that
$$
\|u(0)-v(0)\|_E=0.
$$
So, reasoning as in Theorem \ref{main}, we obtain a unique mild solution of \eqref{abstract-intro} under $(f5\infty)$. Extending such solution as above, gets a unique solution of problem \eqref{Abstract2} associated to \eqref{fixed-initial-condition-abstract}.
\end{proof}

\noindent We are now able to prove the existence of a periodic solution for problem \eqref{periodic-infty}. Consider problem \eqref{periodic-infty} under the following assumptions:
\begin{itemize}
\item[$(h1\infty)$] for every $v \in \mathbb{R}$, $h(\cdot,\cdot,v)\colon \R_+ \times \Omega \to \mathbb{R}$ is a measurable function;
\item[$(h2\infty)$] for every $t \in \R_+$ and $x \in \Omega$, $h(t,x,\cdot) \colon \mathbb{R} \to \mathbb{R}$ is continuous;
\item[$(h3\infty)$] there exists $m > 0$ and $\ell: \mathbb{R}_+ \times \Omega \to \mathbb{R}_+$ such that
\begin{itemize}
    \item[-] $\ell(\cdot,x) \in L^\infty(\mathbb{R}_+;\mathbb{R}_+)$ for a.e. $x \in \Omega$;
    \item[-] $\ell(t,\cdot) \in L^q(\Omega;\mathbb{R}_+)$ for a.e. $t \in \mathbb{R}_+$
    \end{itemize} and such that
$$ 
|h(t,x,v)| \leq \ell(t,x) + m |v|^{p/q}, \quad \mbox{for every} \quad (t,x,v) \in \R_+ \times \Omega \times \mathbb{R},
$$
with $2 \leq q < p < \infty$ for $k \leq 2$ and $2 \leq q < p < \infty$, $\displaystyle\frac{pq}{p-q} > \frac{k}{2}$, for $k > 2$;
\item[$(h4\infty)$] $v \; h(t,x,v) \leq 0$ for a.e. $t \in [0,T]$, $x \in \Omega$ and for every $v \in \mathbb{R}$;
\item[$(h5\infty)$] for every $x \in \Omega$ and $v \in \mathbb{R}$ the map $h(\cdot,x,v):\R_+ \to \R$ is $T$-periodic;
\item[$(h6\infty)$] $(u-v) \; (h(t,x,u)-h(t,x,v)) \leq 0$ for a.e. $(t,x) \in [0,T] \times \Omega$ and for every $u,v \in \mathbb{R}$;

\end{itemize}

\begin{theorem}
\label{ex-pde-infty}
Under the assumptions $(h1\infty)$-$(h5\infty)$, the problem \eqref{periodic-infty} admits at least one solution $u \in C(\R_+;L^p(\Omega; \mathbb{R}))$ such that $\|u(t)\|_p < R_0$ for every $t\ge 0$ for a suitable $R_0 > 0$. Moreover, if $(h6\infty)$ is satisfied we have the uniqueness of mild solutions for the problem \eqref{periodic-infty} with the restriction 
\begin{equation*}
u(0,x)=u_0(x)=u(T,x) \quad \mbox{for a.e.} \quad x \in \Omega
\end{equation*}
where $u_0 \in L^p(\Omega;\mathbb{R})$.
\end{theorem}
\begin{proof}
Following the same reasonings as in Section \ref{existence-pde} we can write problem \eqref{periodic-infty} as the abstract problem \eqref{Abstract2} with $f:\R_+ \times L^p(\Omega; \mathbb{R}) \to L^q(\Omega; \mathbb{R})$ the Nemytskii operator associated to $h$. As in the proof of Theorem \ref{ex-pde1} it is possible to prove that $f$ satisfies the assumptions $(f1\infty)$,$(f3\infty)$, $(f4\infty)$ and $(f5\infty)$. Moreover, condition $(h5\infty)$ trivially implies the periodicity of $f$ with respect to the first argument. Thus we get the existence and uniqueness of periodic mild solution applying Theorem \ref{main2}.
\end{proof}


\noindent {\bf Acknowledgements:} We would like to thank the anonymous referee, who pointed out the question of the uniqueness, thus improving the quality of the paper.

\noindent The research is carried out within the national group GNAMPA of INdAM.

\noindent The first author is partially supported by the projects ``Metodi della Teoria dell'Approssimazione, Analisi Reale, Analisi Nonlineare e loro applicazioni" and ``Integrazione, Approssimazione, Analisi Nonlineare e loro Applicazioni", funded by the 2018 and 2019 basic research fund of the University of Perugia and by a 2020 GNAMPA-INDAM Project ``Processi evolutivi con memoria descrivibili tramite equazioni integro-differenziali". 

\end{document}